\documentclass[11pt,a4paper]{article}
\usepackage{amsmath,amsthm, amssymb, url, enumerate}
\usepackage{graphicx}
\usepackage{cases}
\usepackage{stmaryrd}
\usepackage[T1]{fontenc}
\usepackage{textcomp}
\usepackage{color}
\usepackage[all]{xy}
\usepackage{cases}
\usepackage{stmaryrd}
\usepackage{bm}
\usepackage{here}

\theoremstyle{plain}
\newtheorem{defi}{Definition}[section]
\newtheorem{thm}[defi]{Theorem}
 
\newtheorem{lem}[defi]{Lemma}
\newtheorem{rem}[defi]{Remark}
\newtheorem{cor}[defi]{Corollary}

\newtheorem{ex}[defi]{Example}

\renewcommand{\proofname}{Proof.\ }

\newcommand{\slmc}[1]{{\rm SL}(#1,\mathbb{C})}

\newcommand{\compl}[1]{\mathbb{C}^{#1}}
\newcommand{\real}[1]{\mathbb{R}^{#1}}
\newcommand{\ratio}[1]{\mathbb{Q}^{#1}}
\newcommand{\inte}[1]{\mathbb{Z}^{#1}}
\newcommand{\natu}[1]{\mathbb{N}^{#1}}

\newcommand{\age}{\mathrm{age}}

\newcommand{\CFF}{\mathrm{CFF}}

\makeatletter
\def\mapstofill@{%
   \arrowfill@{\mapstochar\relbar}\relbar\rightarrow}
\newcommand*\xmapsto[2][]{%
   \ext@arrow 0395\mapstofill@{#1}{#2}}
\makeatother

\begin{document}
\title{\vspace{-3cm}Existence of Crepant resolutions for two-parameter\\ Gorenstein Cyclic Quotient Singularities}
\author{Yusuke\ Sato}
\date{}
\maketitle 
\tableofcontents
\begin{abstract}
In this paper, we show a condition for two-parameter Gorenstein cyclic quotient singularities to have a crepant resolution by using the remainder polynomial in any dimension.
\end{abstract}

\section{Introduction}
\noindent
Let $G$ be a finite subgroups of $\slmc{n}$, then $\compl{n}/G$ has a Gorenstein quotient singularity. In the case of $n=2$ and $3$, it is known that all Gorenstein quotient singularities possess crepant resolutions. In the case of $n\geq4$, Gorenstein quotient singularities do not necessarily have a crepant resolution. On the other hand,  D. I. Dais, U.U. Haus and M.Henk have proposed a condition for $\compl{n}/A$ where $A=\frac{1}{r}(a,b,1,\dots,1)$ with $r=a+b+(n-2)$ to possess a  crepant resolution for all dimension \cite{DHH}. We call this two-parameter Gorenstein cyclic quotient singularities.
After that, a new criterion for these quotient singularities to admit a crepant resolution is introduced by S. Davis , T. Logvinenko and M. Reid.
In this paper, we will give the remainder polynomial version of their result.\\

\noindent
\textbf{Theorem \ref{main}.}\ 
{\it Let $\compl n/A$ be a Gorenstein cyclic quotient singularity of type $A=\frac{1}{r}(1,d,c,\dots,c)$.
$\compl n/A$ has a crepant resolution if and only if the ages of all coefficients of the remainder polynomial $\mathcal{R}_*\left( \frac{(1,d,c,\dots,c)}{r} \right)$ are $1$.}\\

The remainder polynomial is introduced by T. Ashikaga as a dimensional generalization of continuous fraction. This multi-dimensional continued fraction consist of the {\it remainder polynomial} and the {\it round down polynomial}, and the round down polynomial is a dimensional generalization of the Hirzebruch-Jung continuous fraction. The remainder polynomial indicates the types of quotient singularities which appear in each step of {\it Fujiki-Oka resolution} ( see \cite{Ashikaga} for more details). 
K. Sato and the author showed the Fujiki-Oka resolution is crepant if and only if the the ages of all coefficients of the remainder polynomial are $1$ \cite{Sato}. 
The result of this paper is derived from the above property of the remainder polynomial.\\
In section $2$, we introduce some notation of toric geometry and the definition of a remainder polynomial. In section $3$, our main results are stated and proved.

\section{Toric geometry and Continued fractions}
First, we set up notion and terminology of toric geometry. 
Let $G$ be a finite cyclic subgroup of ${\rm SL}(n,\compl{})$ generated by ${\rm diag}(\varepsilon^{a_1},\varepsilon^{a_2},\dots,\varepsilon^{a_n})$ where $\varepsilon$ is the $r$-th root of unity. For simplicity, ${\rm diag}(\varepsilon^{a_1},\varepsilon^{a_2},\dots,\varepsilon^{a_n})$ is denote by $\frac{1}{r}(a_1,a_2,\dots,a_n)$.\\

\subsection{Fundamentals on Toric Geometry}\label{TG}
Let $N$ be a free $\inte{}$-module of rank $n$ and $N_\real{} = N \otimes_{\inte{}}\real{}$. Let $\bm{e}_1,\ldots, \bm{e}_n$ be the canonical basis of the vector space $N_\real{}$. We set a rational strongly convex polyhedral cone $\sigma$ as $\real{}_{\geq0}\bm{e}_1 + \cdots + \real{}_{\geq0}\bm{e}_n$ where $\real{}_{\geq0}$ is the set of all non negative elements in $\real{}$. To shorten notation, $\sigma$ also signifies the finite fan consists of all faces of $\sigma$. The {\it dimension} of a cone $\sigma$ is defined as the dimension of $\real{}\cdot\sigma$ as vector space over $\real{}$.  The toric variety $X(N,\sigma)$ determined by $N$ and the finite fan $\sigma$ is isomorphic to $\mathbb{C}^n$. There exists a morphism of toric varieties $\phi_T : X(N,\sigma) \to X(N^{\prime},\sigma)$  corresponding to the quotient map $\phi: \mathbb{C}^n \to \mathbb{C}^n/G$ where $N^{\prime}$ is the free $\inte{}$-module of rank $n$ satisfying $N^{\prime}/N \cong G$ as groups. Therefore, there is a primitive element $\bar{g}=\frac{1}{r}(a_1,\dots,a_n) \in N^{\prime}$ for every $g\in G$.

\begin{defi}\upshape
We define the {\it age} of an element $b=(b_1,\dots,b_n)$ of $N^{\prime}$ to be 
$$
\age(b)=\sum_{i=1}^{n}b_i.
$$
\end{defi}

We shall recall the definition of crepant resolutions as toric geometry.
If a fan $\Sigma$ subdivides the fan $\sigma$, then we have a birational map $f: X(N^{\prime},\Sigma) \to X(N^{\prime},\sigma)$, and the following relation holds between the canonical divisors:
$$
K_{X(N^{\prime},\Sigma)}=f^{*}(K_{X(N^{\prime},\sigma)})+\sum_{\tau \in \Sigma(1)}a_{\tau}D_{\tau},
$$
where $D_{\tau}$ is an exceptional divisor corresponding to the one dimensional cone $\tau \in \Sigma(1)$ in $\Sigma$ and $a_{\tau}=\age(A_{\tau})-1$, where $A_{\tau}$ is the primitive element in $\tau$. The rational number $a_{\tau}$ is called the {\it discrepancy} of $D_{\tau}$. 

\begin{rem}\label{rem2.2}\upshape
Let $\Sigma$ is a subdivision of $\sigma$ by using lattice points of which ages are $1$. If the toric variety $X(N^{\prime},\Sigma)$ is smooth, then $X(N^{\prime},\Sigma)$ is a crepant resolution of $\compl n/G$.
\end{rem}

The convex hull $\mathfrak{s}_G\subset N'_{\real{}}$ spanned by ${\bm e}_1,{\bm e}_2,\ldots ,{\bm e}_n$ is called the {\it junior simplex}. By Remark \ref{rem2.2}, a crepant resolution $X(N^{\prime},\Sigma)$ can be identified with a basic triangulation of  $\mathfrak{s}_G$ by using points in $N'$.\\ 
In this paper, we consider two-parameter cyclic quotient singularities $\compl{n}/A$ where $A$ denote a cyclic group generated by $\frac{1}{r}(a,b,1,\dots,1)$.
These singularities have the following three cases.
\begin{itemize}
\item[(1)] ${\rm GCD}(r,a,b)=d > 1$
\item[(2)] ${\rm GCD}(r,a,b)=1$,  ${\rm GCD}(r,a)=d_1 > 1$ or ${\rm GCD}(r,b)=d_2 > 1$ 
\item[(3)] ${\rm GCD}(r,a)=1$ or ${\rm GCD}(r,b)=1$
\end{itemize}
If $A$ satisfies (1), it is easily seen that $\compl{n}/A$ has a crepant resolution (see \cite{DLR}).
In the case of (2), $\compl{n}/A$ has a crepant resolution if and only if lattice points $\frac{1}{r}(0,k_1,r_1,\dots,r_1)$ and $\frac{1}{r}(k_2,0,r_2.\dots,r_2)$ are on the junior simplex with $r=r_i\cdot d_i$ and $r=k_i+r_i(n-2)$ for $i=1,2$.\\
From now on, we assume that ${\rm GCD} (r,a)=1$. In other words, we treat only the case $A=\frac{1}{r}(1,d,c,\dots,c)$ with $r=1+d+(n-2)c$.


\subsection{The remainder polynomial}\label{ACF}
In this subsection, we shall introduce the definition of the remainder polynomial. 

\begin{defi}\upshape\label{profrac}
Let $n$ be an integer greater than or equal to $1$. Let $\mathbf{a}=(a_1,\dots,a_n) \in \inte n$ and $r \in \natu{}$ which satisfies $0\leq  a_i \leq r-1$ for $1\leq i \leq n$. We call the symbol
$$
\frac{\mathbf{a}}{r}=\frac{(a_1,\dots,a_n)}{r}
$$
an {\it $n$-dimensional proper fraction}.
\end{defi}

\begin{defi}\upshape 
Define the {\it age} of an $n$-dimensional proper fraction $\frac{\mathbf{a}}{r}=\frac{(a_1,\dots,a_n)}{r}$ to be
$$
\age\left( \frac{\mathbf{a}}{r} \right)=\frac{1}{r}\sum_{i=1}^n a_i.
$$
\end{defi}

In the following, the symbol $\ratio{prop}_n$ (resp. $\overline{\ratio{prop}_n}$) means the set of $n$-dimensional proper fractions (resp. the set $\ratio{prop}_n \cup \{\infty\}$). Moreover, $\overline{\ratio{prop}_n}[x_1,\dots,x_n]$  denotes the set consists of all noncommutative polynomials with $n$ variables over $\overline{\ratio{prop}_n}$. 
The remainder polynomials is obtained via  {\it remainder maps} for a {\it semi-unimodular proper fraction} (i.e., a proper fraction such that at least one component of $\mathbf{a}$ is $1$). Roughly speaking, the remainder map is division for just one component of the vector $\mathbf{a}$ by $r$.

\begin{defi}\upshape (\cite[Def 3.1.]{Ashikaga})
Let $\frac{(a_1,\dots,a_n)}{r}$ be a semi-unimodular proper fraction. 
For $1\leq i \leq n$, the {\it $i$-th remainder map} $R_i:\overline{\ratio{prop}_n} \to \overline{\ratio{prop}_n}$ is define by
$$
R_i\left(\frac{(a_1,\dots,a_n)}{r}\right)=
\left\{ 
\begin{array}{cc}
 \left(\frac{\overline{a_1}^{a_i},\ \dots,\ \overline{a_{i-1}}^{a_i},\ \overline{-r}^{a_i},\ \overline{a_{i+1}}^{a_i},\ \dots, \overline{a_n}^{a_i}}{a_i}
         \right)& {\rm if} \  a_i \neq 0\\  
 \infty &  {\rm if}\  a_i=0
\end{array}
\right.
$$
and $R_i(\infty)=\infty$ where $\overline{a_j}^{a_i}$ is an integer satisfying $0\leq \overline{a_j}^{a_i} < a_i$ and $\overline{a_j}^{a_i} \equiv a_j$ modulo $a_i$.
\end{defi}

\begin{ex}\upshape
If $v=\frac{(1,2,5,7)}{8}$, then 
\begin{eqnarray*}
R_2(v)&=&\frac{(1,0,1,1)}{2}\ \text{and}\\
R_3(v)&=&\frac{(1,2,2,2)}{5}.
\end{eqnarray*}

\end{ex}

\begin{defi}\upshape\cite[Def 3.2.]{Ashikaga}\label{DOACF}
Let $\frac{\mathbf{a}}{r}$ be an $n$-dimensional semi-unimodular proper fraction,  and $\mathbf{I}=\{1,2,\dots,n\}$ signifies the index set of the variables. \\
The {\it remainder polynomial} $\mathcal{R}_*\left(\frac{\mathbf{a}}{r}\right) \in \overline{\ratio{prop}_n}[x_1,\dots,x_n]$ is defined by
  $$
  \mathcal{R}_*\left(\frac{\mathbf{a}}{r}\right)=\frac{\mathbf{a}}{r}+
                                      \sum_{(i_1,i_2,\dots,i_l)\in \mathbf{I}^l,\: l\geq 1 }(R_{i_l}\cdots R_{i_2}R_{i_1})\left(\frac{\mathbf{a}}{r}\right)\cdot x_{i_1}x_{i_2}\cdots x_{i_l}
  $$
  where we exclude terms with coefficients $\infty$ or $\frac{(0,0,\dots,0)}{1}$.
\end{defi}

\begin{defi}\upshape
The term with the variable $x_i\cdots x_i$ in a remainder polynomial is called to be {\it iterated} where $1 \leq i \leq n$, and 
the lattice point in $N'$ corresponding to the coefficient of iterated terms is also called to be {\it iterated}.
\end{defi}

\begin{ex}\upshape
Let $v=\frac{(1,2,6,6)}{15}$, then the remainder polynomial is
\begin{eqnarray*}
\mathcal{R}_*\left(\frac{(1,2,6,6)}{15}\right)&=&
             \frac{1}{15}(1,2,6,6)+\frac{1}{2}(1,1,0,0)x_2+\frac{1}{6}(1,2,3,0)x_3+\frac{1}{6}(1,2,0,3)x_4 \\
                                &+&\frac{1}{2}(1,0,1,0)x_3x_2+\frac{1}{3}(1,2,0,0)x_3x_3+\frac{1}{2}(1,0,0,1)x_4x_2+\frac{1}{3}(1,2,0,0)x_4x_4\\
                                &+&\frac{1}{2}(1,1,0,0)x_3x_3x_2+\frac{1}{2}(1,1,0,0)x_4x_4x_2
\end{eqnarray*}
Moreover, iterated terms are $\frac{1}{15}(1,2,6,6),\frac{1}{2}(1,1,0,0)x_2,\frac{1}{6}(1,2,3,0)x_3,\frac{1}{6}(1,2,0,3)x_4, \frac{1}{3}(1,2,0,0)x_3x_3$ and $\frac{1}{3}(1,2,0,0)x_4x_4$ .
\end{ex}

\section{The main result}
\subsection{The criterion in the remainder polynomial form}
In this section, we shall show a condition for two-parameter Gorenstein cyclic quotient singularities to have a crepant resolution by using the remainder polynomial. To prove the main result, we introduce some properties of the remainder polynomial without proof. These properties was proved by the author and K. Sato(\cite{Sato}, Section 3).

\begin{lem}\label{Lem1}\upshape (\cite[Lemma 3.4.]{Sato})
Assume that $1+a_2+a_3+\dots+a_n=r$ for $\frac{1}{r}(1,a_2,\dots,a_n)$. Then $\age\left(\mathcal{R}_i\left(\frac{(1,a_2,\dots,a_n)}{r}\right)\right)$ is an integer.
\end{lem}

\begin{thm}\upshape (\label{Satothm}\cite[Theorem 3.1., Proposition 3.9.]{Sato})
Let $\compl n/G$ be a quotient singularity of $\frac{1}{r}(1,a_2,\dots,a_n)$-type satisfying $1+a_2+\cdots+a_n=r$. 
\begin{itemize}
\item[(i)] If the ages of all coefficients of $\mathcal{R}_*\left( \frac{(1,a_2,\dots,a_n)}{r} \right)$ are $1$, then the Fujiki-Oka resolution $X(N^{\prime},\CFF(\sigma))$ gives a crepant resolution for $X(N^{\prime},\sigma)$.
\item[(ii)] If the remainder polynomial $\mathcal{R}_*\left(\frac{(1,a_2,\dots,a_n)}{r}\right)$ contains an iterated term of which the age of the coefficient is equal to or bigger than $2$, then $\compl n/G$ has no crepant resolutions.
\end{itemize}
\end{thm}

\begin{cor}\label{maincor}\upshape (\cite[Corollary 3.5.]{Sato}, )
For all three dimensional semi-isolated Gorenstein quotient singularities, the Fujiki-Oka resolutions are crepant.
\end{cor}

\proofname Let $G=\left\langle \frac{1}{r}(1,a,b)\right \rangle$ where $1+a+b=r$. we have $\mathcal{R}_2\left(\frac{(1,a,b)}{r}\right)=\frac{(1,\overline{-r}^a,\overline{b}^a)}{a}$, and the age of $\mathcal{R}_2\left(\frac{(1,a,b)}{r}\right)$ is an integer by Lemma \ref{Lem1}.
Clearly, $1+\overline{-r}^a+\overline{b}^a< 2a $. So, the age of $\mathcal{R}_2\left(\frac{(1,a,b)}{r}\right)$ equals to $1$.
Thus, the ages of all coefficients of $\mathcal{R}_*\left(\frac{(1,a,b)}{r}\right)$ equal to $1$. By Theorem \ref{Satothm}, the Fujiki-Oka resolution $X(N^{\prime},\CFF(\sigma))$ is crepant.
\qed\\

Applying the Theorem \ref{Satothm} to the cyclic group $A=\frac{1}{r}(1,d,c,\dots,c)$ gives the conditions to the existence of crepant resolutions.

\begin{lem}\label{keylem}\upshape
If  {\large $\mathcal{R}_*\left(\frac{(1,d,c,\dots,c)}{r}\right)$} with $1+d+(n-2)c=r$ does not satisfy the condition (ii), then {\large $\mathcal{R}_*\left(\frac{(1,d,c,\dots,c)}{r}\right)$} satisfies the condition (i).
\end{lem}

\proofname
It is easily to check that the age of $R_i\left(\frac{(1,d,c,\dots,c)}{r}\right)=\frac{(1,\overline{d}^c,0,\dots, 0,\overline{-r}^c,0,\dots,0)}{c}$ is equal to $1$ for $i=3,\dots,n$. By the proof of Corollary\ref{maincor}, 
$\mathcal{R}_*\left(R_i\left(\frac{(1,d,c,\dots,c)}{r}\right)\right)$
satisfies the condition (i). On the other hand, by assumption, the image of the remainder map $(R_{2}\cdots R_{2})\left(\frac{(1,d,c,\dots,c)}{r}\right)$ is $\frac{1}{r'}(1,d',c',\dots,c')$ for some positive integer $r', d', c'$ with $1+d'+(n-2)c'=r'$.
and  $\mathcal{R}_*\left(R_i\left(\frac{(1,d,c,\dots,c)}{r}\right)\right)$, Thus, $\mathcal{R}_*\left(R_i\left(\frac{(1,d',c',\dots,c')}{r'}\right)\right)$ satisfies the condition(i) for $i=3,\dots,n$. By induction, it follows that $\mathcal{R}_*(\frac{(1,d,c,\dots,c)}{r})$ satisfies the condition (i).  
\qed

Lemma \ref{keylem} and Theorem \ref{Satothm} lead to the following theorem.

\begin{thm}\label{main}
Let $\compl n/G$ be a quotient singularity of $\frac{1}{r}(1,d,c,\dots,c)$-type.
$\compl n/G$ has a crepant resolution if and only if the ages of all coefficients of the remainder polynomial $\mathcal{R}_*\left( \frac{(1,d,c,\dots,c)}{r} \right)$ are $1$.
\end{thm}

\subsection{Relationship with Hirzebruch-Jung continued fraciton}
At the end of this paper, we discuss the relationship between the above results and the following condition given by S. Davis , T. Logvinenko and M. Reid \cite{DLR}:

\begin{thm}\label{DLR}\upshape (\cite[Theorem 2.5]{DLR})
Let $A=\left\langle \frac{1}{r}(1,d,c,\dots,c)\right\rangle \subset \slmc{n}$ with $1+d+(n-2)c=r$ and ${\rm GCD}(r,d)=1$, then the following two conditions are equivalent.
\begin{itemize}
\item[(a)] There exists a crepant resolution.
\item[(b)] The Hirzebruch-Jung expansion of $\frac{r}{d}$ has every entry congruent to $2$ modulo $n-2$.
\end{itemize}
\end{thm}
By the proof of Lemma\ref{keylem}, it is enough to check existence of crepant resolutions that we only calculate the iterated terms $(R_{2}\cdots R_{2})\left(\frac{(1,d,c,\dots,c)}{r}\right)\cdot x_{2}^i$ of the remainder polynomial.
In other words, we have the following condition (c) and the conditions (a), (b) and (c) are equivalent.
\begin{itemize}
\item[(c)] The age of $R^i\left(\frac{(1,d,c,\dots,c)}{r}\right)$
equals to $1$ for all $i$, where $R^i\left(\frac{(1,d,c,\dots,c)}{r}\right)$ denote the coefficient of the iterated term $(R_{2}\cdots R_{2})\left(\frac{(1,d,c,\dots,c)}{r}\right)\cdot x_{2}^i$ with $i=1,\dots,s$.
\end{itemize}
 
We introduce directly the proof (c) $\Rightarrow$ (b). Assume that the Hirzebruch-Jung expansion of $\frac{r}{d}$ has entry $a_1,\dots,a_s$.
It is sufficient to show that $a_1$ is congruent to $2$ modulo $n-2$ when the condition (c) holds for $i=1$. 
For $\frac{(1,d,c,\dots,c)}{r}$, the image of the $2$nd remainder map is 
$$R_2\left(\frac{(1,d,c,\dots,c)}{r}\right)=\frac{(1,\overline{-r}^d,\overline{c}^d,\dots,\overline{c}^d)}{d}.$$ 
By assumption, $\age(R_2(\frac{(1,d,c,\dots,c)}{r}))=1$ holds. Thus we have the followings equation.
\begin{eqnarray}
1+\overline{-r}^d+(n-2)\overline{c}^d=d \\
1+d+(n-2)c=r
\end{eqnarray}
Subtracting the equation (1) from (2), we obtain
\begin{eqnarray}
(c-\overline{c}^d)(n-2)+2d=r+\overline{-r}^d=-d \cdot q
\end{eqnarray}
 where $q$ denote the quotient of $-r$ divided by $d$. Since $(c-\overline{c}^d)$ is divisible by $d$, we have $2\equiv -q \quad ({\rm mod}\  n-2)$. By definition of the quotient, $\frac{-r}{d}=q+\frac{\overline{-r}^d}{d}$ holds. This leads to $-q=a_1$.
\qed



\textsc{Graduate School Of Mathematical Sciences, University Of Tokyo 3-8-1 Komaba, Meguro-ku, Tokyo 153-8914, Japan.}\\
E-mail:yusuke.sato@ipmu.jp


\begin{thebibliography}{}
{\small
\bibitem[A]{Ashikaga} T. Ashikaga, {\it Multidimensional continued fractions for cyclic quotient singularities and Dedekind sums}, To appear in Kyoto J. Math. Advance publication (2019).

\bibitem[DHH]{DHH} D. I. Dais, U.U. Haus and M.Henk, {\it On crepant resolutions of $2$-parameter series of Gorenstein cyclic quotient singularities}, Results Math., {\bf 33} (1998) 208-265.
\bibitem[DLR]{DLR} S. Davis, T. Logvinenko, M. Reid, How to calculate $A{\rm -Hilb}(\mathbb{C}^n)$ for $\frac{1}{r}(a,b,1,\dots,1)$, preprint.
\bibitem[F]{Fujiki} A. Fujiki, {\it On resolution of cyclic quotient singularities}, Publ. Res. Inst. Math. Sci. {\bf 10} (1974/75) 293--328.
\bibitem[O]{Oka} M. Oka, {\it On the resolution of hypersurface singularities}, Adv. St. Pure Math. {\bf 8} (1986) 405--436.
\bibitem[SS]{Sato} K. Sato and Y. Sato, {\it Crepant Property of Fujiki-Oka Resolutions for Gorenstein Abelian Quotient Singularities}, preprint, arXiv:math/2004.03522
}
\end{thebibliography}
\end{document}